\documentclass[a4paper,12pt]{article}
\usepackage{amsmath,amssymb,amsthm,amsfonts,epsfig}

\numberwithin{equation}{section}

\newtheorem{thm}{Theorem}[section]

\newtheorem{rem}{Remark}[section]
\newtheorem{prop}{Proposition}[section]
\newtheorem{cor}{Corollary}[section]

\title{A priori estimate for a family of semi-linear elliptic equations
with critical nonlinearity}

\author{
{\sc Lei Zhang} \\
 {\sc \footnotesize  University of Alabama at Birmingham}  \\
 {\sc \footnotesize
    Department of Mathematics,}\\
    {\sc \footnotesize 452 Campbell Hall} \\
 {\sc \footnotesize
    1300 University Boulevard}\\
 {\sc \footnotesize
    Birmingham, AL 35294-1170}  \\
{\tt \footnotesize
    leizhang@math.uab.edu}\thanks{Supported by
National Science Foundation Grant 0600275 (0810902).} }

\date{}
\input { amssym.def}

\begin{document}
\maketitle

\begin{abstract}
We consider positive solutions of $\Delta u-\mu
u+Ku^{\frac{n+2}{n-2}}=0$ on $B_1$ ($n\ge 5$) where $\mu $ and
$K>0$ are smooth functions on $B_1$. If $K$ is very sub-harmonic
at each critical point of $K$ in $B_{2/3}$ and the maximum of $u$
in $\bar B_{1/3}$ is comparable to its maximum over $\bar B_1$,
then all positive solutions are uniformly bounded on $\bar
B_{1/3}$. As an application, a priori estimate for solutions of
equations defined on $\Bbb S^n$ is derived.
\end{abstract}

\noindent


 \smallskip

  {\bf Running title:} A priori estimate for semi-linear
  equations.


\section{Introduction}

In this article we study the equation
\begin{equation}
\Delta u-\mu(x) u+K(x)u^{n^*}=0,\quad u>0,\quad u\in C^2(B_1),
\quad n\ge 5, \label{feb14e1}
\end{equation}
where $B_1$ is the unit ball centered at the origin,
$n^*=\frac{n+2}{n-2}$ is the critical power in Sobolev embedding,
$\mu$ is a $C^1$ function on $B_1$ and $K\in C^3(B_1)$ is a
positive function. We shall derive a priori estimate under natural
assumptions on $K$ and $\mu$.

Equation (\ref{feb14e1}) has rich connections in physics and
geometry. In particular, it is very closely related to the well
known Yamabe equation, which has been extensively studied for
decades. Many interesting features of the Yamabe equation are also
reflected on this equation. When $\mu$ is a non positive constant,
many results on the existence of solutions, multiplicity of
solutions, a priori estimates, bifurcation phenomena, Harnack type
inequalities, etc can be found in the literature. We refer the
interested readers to
\cite{BG}\cite{BN}\cite{cerqueti}\cite{Ouyang}\cite{zhu} and the
references therein. On the other hand, much less references can be
found for the case $\mu>0$.
 A recent paper of Lin-Prapapet \cite{linprajapat}
discussed the case $\mu=constant>0$ and they pointed out that it
is also interesting and important to study the following Harnack
type inequality:
\begin{equation}
\label{sept2e1}
 (\max_{\bar B_{1/3}}u)( \min_{\bar B_{2/3}} u)\le C.
 \end{equation}
In a slightly different setting, they derived this Harnack
inquality for $3\le n\le 6$ under some flatness assumptions of
$\nabla K$ near its critical points. They also speculated that
(\ref{sept2e1}) should still hold for higher dimensions under
similar assumptions by finer analysis.

The Harnack inequality (\ref{sept2e1}) is an important estimate to
understand the blowup phenomena of (\ref{feb14e1}) and many
related equations with the critical Sobolev nonlinearity. The very
first discussions of this inequality can be found in \cite{S1} and
\cite{CSL1}. With this Harnack type inequality, usually the blowup
phenomena is greatly simplified and some energy estimates are
implied. Moreover, some further results such as a priori estimate,
precise description of the blowup bubbles, etc can be obtained.

In this article, we use a very different approach from
Lin-Prajapat's to obtain a priori estimate for general
 $\mu\in C^1$ and $n\ge 5$. We shall also derive the Harnack type inequality
 as an intermediate step toward our result. Our idea stems from the author's
joint work with Y. Y. Li \cite{LiZhang1}\cite{LiZhang3} on the
compactness of solutions of the Yamabe equation.

We assume the following on $K$ and $\mu$:
$$(K,\mu):\quad C_1^{-1}\le K(x)\le C_1, \quad
x\in B_1,\quad \|K\|_{C^3(B_1)}\le C_1,\quad \|\mu
\|_{C^1(B_1)}\le C_0. $$

In addition, we need the maximum of $u$ in $B_{1/3}$ comparable to
its maximum in $B_1$: There exists $C_2>0$ such that
\begin{equation}
\label{aug15e1} \max_{\bar B_{1/3}}u\ge \frac 1{C_2}\max_{\bar
B_{1}}u.
\end{equation}

The main result of the paper is

\begin{thm}
\label{thm1} Given $(K,\mu )$ and (\ref{aug15e1}) there exists
$C_3(C_0, n, C_1, C_2)>0$ such that for each critical point $x$ of
$K$ in $B_{\frac 23}$, if $\Delta K(x)>C_3$, any solution $u$ of
(\ref{feb14e1}) satisfies
\begin{equation}
\label{feb14e2}
 \max_{\bar B_{1/3}}u\le C_4(C_0,n,C_1,C_2).
\end{equation}
\end{thm}

If $\mu $ is a positive constant, Lin-Prajapet proved in
\cite{linprajapat} without the assumption (\ref{aug15e1})
 that for $n=3$, if $K$ is
H\"older continuous with exponent $\theta\in (\frac 12, 1]$ then
(\ref{sept2e1})
 holds. For $4\le n\le 6$, if $K\in C^1$ and in a
neighborhood of each critical point $x_0$ of $K$,
$$c|x-x_0|^{\theta-1}\le |\nabla K(x)|\le C|x-x_0|^{\theta-1}$$
holds for $c, C>0$ and $\frac{n-2}2\le \theta\le n-2$, then
(\ref{sept2e1}) holds.

Even though (\ref{aug15e1}) seems to be a strong restriction for
locally defined equations, it can be removed when Theorem
\ref{thm1} is applied to globally defined equations. In this case
Theorem \ref{thm1} is particularly useful. For example, let $w$ be
a positive solution to
\begin{equation}
\label{8feb25e1} Lw-aw+Rw^{\frac{n+2}{n-2}}=0,\quad \mbox{on}\quad
S^n
\end{equation}
where $L=\Delta_{g_0}-\frac{n(n-2)}4$ is the conformal Laplacian
operator of $(S^n,g_0)$, $a$ and $R$ are positive smooth
functions. By using the stereographic projection $\pi$ from $S^n$
to $\Bbb R^n$, we set $K(y)=R(\pi^{-1}(y))$ for $y\in \Bbb R^n$
and $\mu(y)=a(\pi^{-1}(y))$. Without loss of generality we assume
the north pole is not a critical point of $R$. Similar to
$(K,\mu)$ we assume
$$(K,\mu)_1:\quad C_5^{-1}\le K\le C_5,\quad \|K\|_{C^3(\Bbb
R^n)}\le C_5,\quad \|\mu\|_{C^1(\Bbb R^n)}\le C_5.$$ Then we have
\begin{cor}\label{cor1}
Given $n\ge 5$ and $(K,\mu)_1$, then there exists $C_6(C_5,n)>0$
such that for each critical point $y$ of $K$, if $\Delta
K(y)>C_6$, solution $w$ of (\ref{8feb25e1}) satisfies
$$C_7^{-1}\le w(x)\le C_7,\quad x\in \Bbb S^n $$
where $C_7(C_5,n)>0$ is independent of $w$.
\end{cor}

\begin{rem}\label{feb26r1}
If $a\equiv 0$ in (\ref{8feb25e1}), then we only need to assume
$\Delta K(y)>0$ on critical point $y$ to get the same a priori
estimate.
\end{rem}

 The results in Corollary \ref{cor1} and Remark \ref{feb26r1} can be
 compared with closely
related results in \cite{Chang}\cite{CSL0}\cite{Li1},\cite{Li2},
etc and the references therein. Comparing to these results,
Corollary \ref{cor1} gives the a priori estimate under a very
short assumption on $K$.

The conclusion of Theorem \ref{thm1} can also be compared with
related results on the compactness of solutions of the Yamabe
equation. It is proved in \cite{LiZhang1} and \cite{marques} that
for the Yamabe equation, blowup point can not appear at places
where the Weyl tensor is $0$. The role of $\Delta K$ in Theorem
\ref{thm1} is similar to that of the Weyl tensor for the Yamabe
equation. Please also see \cite{Brendle1}\cite{Brendle2}
\cite{Chang}\cite{CL2}\cite{Khuri}\cite{LiZhang2}\cite{Lizhu} for
related discussions.

In another work of the author \cite{z1}, among other things the
following result is essentially proved:  For $\mu=0$ and $n\ge 5$,
suppose (\ref{aug15e1}) holds and $\Delta K(x)>0$ for all critical
point $x$ in $B_{2/3}$, then $\max_{\bar B_{1/3}}u\le C$ for some
$C>0$.

Based on Theorem \ref{thm1} and the results in
\cite{linprajapat}\cite{z1} we propose the following two
questions for the case $\mu=constant>0$:
\begin{enumerate}
\item For $n=4$, under the assumptions of Theorem \ref{thm1}, can
one obtain the Harnack type inequality (\ref{sept2e1}) or even
stronger, the a priori estimate? \item For $n\ge 5$ under the
assumptions $(K,\mu)$ and (\ref{aug15e1}), what is the smallest
$C>0$ so that Theorem \ref{thm1} holds for $\Delta K>C$ at each
critical point in $B_{2/3}$? We suspect that the best constant
only comes from the Pohozaev Identity.
\end{enumerate}

The idea of the proof of Theorem \ref{thm1} is by an iterative use
of the well known method of moving spheres (MMS). Some estimates
established by Chen-Lin \cite{CSL2} are very crucial to our
approach. The reason that we need to apply MMS
 many times is because we need to construct appropriate
test functions for this well known method. The construction of
test functions depends on the estimates of some error terms and is
closely related to the spectrum of the linearized operator of the
equation. At the beginning we only have crude estimates of these
error terms. As a consequence we can only construct test functions
( and apply MMS) on these small domains. However, MMS and
Chen-Lin's estimate lead to better estimates of error terms, which
make it possible to construct test functions on larger domains.
After applying this procedure iteratively we obtain the desired
estimates on the error terms and MMS can be applied on roughly the
whole domain.

More specifically the outline of the proof is as follows. Suppose
there is no uniform bound for a sequence of solutions $u_i$. Scale
$u_i$ appropriately so that the maximum of the re-scaled function
$v_i$ is comparable to $1$. Then these functions are defined on
very large balls,say, $B(0, \frac 1{10}M_i^{\frac 2{n-2}})$, where
$M_i$ is the maximum of $u_i$. The first step is to show that
$v_i$ is comparable to a standard bubble that takes $1$ at $0$
over the range $B(0,M_i^{\frac 2{(n-2)^2}})$. The way to prove
this step is by an easy application of MMS on this range. Then in
step two we show that the difference between $v_i$ and the
standard bubble is of the order $O(M_i^{-\frac 2{n-2}})$.  The
approach for this step is based on Chen-Lin's argument in
\cite{CSL2}. The result in step two helps us to describe some
error terms in a better way so that we can use MMS to prove, in
step three, that $v_i$ is comparable to the standard bubble over
the range $B(0, M_i^{\frac 4{(n-2)^2}})$. In this step we need to
rewrite some major error terms into a product of spherical
harmonics with radial functions. This decomposition allows us to
find test functions of the same form. By using the Pohozaev
identity in step four we obtain that $|\nabla K|$ at the blowup
point must vanish at the order of $O(M_i^{-\frac 2{n-2}})$. Then
in step five we apply the Chen-Lin estimate again to show that the
closeness between $v_i$ and the standard bubble can be improved to
$O(M_i^{-\frac 4{n-2}})$. This is optimal for this closeness. This
new estimate helps us again to describe some error terms in a
better way so that we can apply MMS again in a bigger domain. In
fact in step six we show that in almost the whole domain, $v_i$ is
comparable to the standard bubble. In this step, the largeness of
$\Delta K$ at the blowup point is used. Finally in step seven we
apply the Pohozaev identity over the whole domain, using symmetry
and all previous estimates to get a contradiction.

The organization of this paper is as follows. In section two the
proof of Theorem \ref{thm1} is presented. The different steps of
the proof are contained in different sub-sections. At the end of
section two we use Theorem \ref{thm1} to prove Corollary
\ref{cor1}.

\section{The proof of Theorem \ref{thm1}}

 In the proof of Theorem \ref{thm1} we only consider the
 case $\mu=constant>0$. In our argument, the
 difference between $\mu$ being a general $C^1$ function and a constant only
 produces  minor terms in our estimate. In order not to make
 notations difficult, we leave the general case for the interested
 readers.

 The proof is based on an assumption for contradiction. Suppose a
sequence of functions $u_i$ can be found to satisfy
$$\Delta u_i-\mu_i u_i+K_i(x)u_i^{n^*}=0 \quad B_1$$
with $K_i$ satisfying $(K, \mu)$ and $0\le \mu_i\le C_0$ such that
\begin{equation}
\label{aug14e9}
u_i(\bar x_i)=\max_{\bar B_{1/3}}u_i\to \infty.
\end{equation}
 Then by a standard selection process we have
$x_i\in B(\bar x_i, 1/10)$ such that $x_i$ is a local maximum of
$u_i$ and $M_i=u_i(x_i)$ is comparable to the maximum of $u_i$ in
$B_1$. Moreover
$$v_i(y)=M_i^{-1}u_i(M_i^{-\frac 2{n-2}}y+x_i)$$
tends in $C^2_{loc}(\Bbb R^n)$ norm to $U$ which satisfies
$$\Delta U+\lim_{i\to \infty}K_i(x_i)U^{\frac{n+2}{n-2}}=0, \quad
\Bbb R^n.$$ and $U(0)=1=\max_{\Bbb R^n} U$. By the well known
classification theorem of Caffarelli-Gidas-Spruck we have
$$U(y)=(1+|y|^2)^{-\frac{n-2}2}.$$
Here we have assumed without loss of generality that
$$\lim_{i\to \infty}K_i(x_i)=n(n-2).$$
The standard selection process can be found in quite a few papers,
for example \cite{z1}. Direct computation shows that the equations
for $v_i$ is
 \begin{equation}
 \label{aug3e2}
 \Delta v_i-\mu_i M_i^{-\frac 4{n-2}}v_i+K_i(M_i^{-\frac 2{n-2}}\cdot
 +x_i)v_i^{n^*}=0, \quad \Omega_i
 \end{equation}
 where $\Omega_i:=B(0,\frac 1{10} M_i^{\frac 2{n-2}})$. In the
 sequel unless we state otherwise a constant always depends on
 $n,C_0,C_1,C_2$.

\subsection{ Estimate of $v_i$ over $B(0,\delta M_i^{\frac
2{(n-2)^2}})$}

In this subsection we establish the following estimate:
\begin{prop}
\label{aug3p1}
 For any $\epsilon>0$ there exists $\delta_0(\epsilon)>0$ such
 that for all large $i$
$$\min_{|y|=r} v_i(y)\le (1+\epsilon)r^{2-n},\quad \forall r
\in (0,\delta_0 M_i^{\frac{2}{(n-2)^2}}).$$
\end{prop}

\noindent{\bf Proof of Proposition \ref{aug3p1}:} The proof is by
a contradiction. Suppose there exists $\epsilon_0>0$ and
$r_i=\circ(1) M_i^{\frac{2}{(n-2)^2}}$ so that
\begin{equation}
\label{aug6e1} \min_{|y|=r_i}v_i(y)\ge (1+\epsilon_0)r_i^{2-n}.
\end{equation}

Note that by the convergence of $v_i$ to $U$ we certainly have
$r_i\to \infty$. We shall use the moving sphere argument. Here we
let $\Sigma_{\lambda}=B_{r_i}\setminus \bar B_{\lambda}$. The
boundary condition for $v_i$ on $|y|=r_i$ is (\ref{aug6e1}). Let
 $$v_i^{\lambda}(y)=(\frac{\lambda}{|y|})^{n-2}v_i(y^{\lambda}),\quad
 y^{\lambda}=\frac{\lambda^2y}{|y|^2} $$
 be the Kelvin transformation of $v_i$ with respect to $\partial
 B_{\lambda}$. Note that in this article $\lambda$ is always
 assumed to stay between two positive constants independent of $i$. The
 equation for $v_i^{\lambda}$ is
 \begin{equation}
 \label{aug3e3}
 \Delta v_i^{\lambda}-\mu_i M_i^{-\frac 4{n-2}}(\frac{\lambda}{|y|})^4v_i^{\lambda}
 +K_i(M_i^{-\frac 2{n-2}}y^{\lambda}+x_i)(v_i^{\lambda})^{n^*}=0\quad
 \Sigma_{\lambda}.
 \end{equation}
Set $w_{\lambda}=v_i-v_i^{\lambda}$ in $\Sigma_{\lambda}$. For
 simplicity we omit $i$ in this notation. We shall apply the
 moving sphere argument to $w_{\lambda}$ with a test function.
 The equation for $w_{\lambda}$ is
 \begin{equation}
 \label{aug3e5}
 T_{\lambda}w_{\lambda}=E_{\lambda}\quad \mbox{in}\quad
 \Sigma_{\lambda}.
 \end{equation}
 where
 $$T_{\lambda}:=\Delta-\mu_i M_i^{-\frac 4{n-2}}+n^*K_i(M_i^{-\frac 2{n-2}}y+x_i)\xi_{\lambda}^{\frac
 4{n-2}},$$
  $\xi_{\lambda}$ is obtained from the mean value theorem:
 $$\xi_{\lambda}^{\frac
 4{n-2}}=\int_0^1(tv_i+(1-t)v_i^{\lambda})^{\frac 4{n-2}}dt,$$
\begin{equation}
 \label{aug3e4}
 E_{\lambda}=\mu_i M_i^{-\frac
 4{n-2}}v_i^{\lambda}(1-(\frac{\lambda}{|y|})^4)
 +(K_i(M_i^{-\frac 2{n-2}}y^{\lambda}+x_i)-K_i(M_i^{-\frac
 2{n-2}}y+x_i))(v_i^{\lambda})^{n^*}
 \end{equation}
 is the main error term.

 For the moving sphere argument we shall find two constants $\lambda_0$ and
 $\lambda_1$, both independent of $i$ such that $\lambda_0\in (\frac
 12, 1)$, $\lambda_1\in (1,2)$. We shall only consider $\lambda\in
 [\lambda_0,\lambda_1]$ for the moving sphere method.
 Test function $h_{\lambda}$, which depends on $i$ ($\lambda\in
 [\lambda_0,\lambda_1]$),
 will be constructed to satisfy
 \begin{eqnarray}
 && h_{\lambda}|_{\partial B_{\lambda}}=0 \label{aug14e2}\\
 &&h_{\lambda}=\circ(1)r^{2-n}\quad \mbox{in}\quad
 \Sigma_{\lambda},\quad |\nabla h_{\lambda}|=\circ(1)\quad \mbox{in}\quad
  \Sigma_{\lambda}\cap B_R,\label{aug14e3}\\
  &&
 T_{\lambda}h_{\lambda}+E_{\lambda}\le 0,\quad \mbox{in}\quad
 O_{\lambda}:=\{y\in \Sigma_{\lambda};\quad v_i(y)\le 2v_i^{\lambda}(y)\quad \}
 \label{aug14e4}
 \end{eqnarray}
 where $R$ in (\ref{aug14e3}) is any fixed large constant.

 Once such a test function is constructed the moving sphere
 argument can be applied to get a contradiction to (\ref{aug6e1}).
 In fact first we show that the moving
 sphere process can get started at $\lambda_0$:
\begin{equation}
\label{aug6e2} w_{\lambda_0}+h_{\lambda_0}>0 \quad \mbox{in}\quad
\Sigma_{\lambda_0}.
\end{equation}
  To see this, first we state a property of the standard bubble $U$:
  \begin{equation}
  \label{aug17e1}
U(y)-U^{\lambda}(y)\sim
(1-\frac{\lambda}{|y|})(1-\lambda)|y|^{2-n},\quad |y|>\lambda,
\end{equation}
which implies that for any $\lambda_0<1$ and $|y|>\lambda_0$
$$w_{\lambda_0}(y)>C(|y|-\lambda_0)|y|^{1-n},\quad
\lambda_0<|y|<\bar R$$  for all large $i$ by the convergence of
$w_{\lambda_0}$ to $U-U^{\lambda_0}$ in $C^2_{loc}(\Bbb R^n)$.
Here $\bar R$ is a large fixed number to be determined.
 By (\ref{aug14e2}) (\ref{aug14e3})
 one sees easily that $w_{\lambda_0}+h_{\lambda_0}>0$
 over $B_{\bar R}\cap \Sigma_{\lambda_0}$. For $|y|>\bar R$ we
 observe that
 $$v_i^{\lambda_0}(y)\le \lambda_0^{n-2}(1+\epsilon_1)|y|^{2-n},
 \quad |y|\ge \bar R$$
 where $\epsilon_1$ is sufficiently small so that
 $\lambda_0^{n-2}(1+\epsilon_1)<1-5\epsilon_1$.
 On the other hand, by the convergence of $v_i$ to $U$ we
 can make $\bar R$ large enough so that
 $$v_i(y)\ge (1-2\epsilon_1)|y|^{2-n},\quad |y|=\bar R.$$
 On $\Sigma_{\lambda_0}\setminus \bar B_{\bar R}$ $v_i$
 satisfies
 $$\Delta v_i-\mu_i M_i^{-\frac 4{n-2}}v_i\le 0, \quad
\Sigma_{\lambda_0}\setminus \bar B_{\bar R}.$$ Let
$$\tilde O:=\{\,\, y\in \Sigma_{\lambda_0}\setminus \bar B_{\bar
R};\quad v_i(y)\le 2|y|^{2-n}\,\, \}. $$

Clearly $w_{\lambda_0}+h_{\lambda_0}>0$ in
$\Sigma_{\lambda_0}\setminus (B_{\bar R}\cup \tilde O)$. Let $G$
be the Green's function of $-\Delta$ on
$\Sigma_{\lambda_0}\setminus \bar B_{\bar R}$ with respect to the
Dirichlet boundary condition, let
$$\phi(y)=\int_{\tilde O}G(y,\eta)\mu_i M_i^{-\frac 4{n-2}}v_i(\eta
)d\eta.$$ Then
$$-\Delta \phi=\mu_i M_i^{-\frac 4{n-2}}v_i \quad \tilde O.$$
Elementary estimate gives
\begin{equation}
\label{aug21e1} \phi(y)\le C(n)\mu_i M_i^{-\frac 4{n-2}}|y|^{4-n}.
\end{equation}
So
\begin{equation}
\label{aug21e2} \phi(y)\le \epsilon_1|y|^{2-n}.
\end{equation}
Now we have
$$\Delta (v_i+\phi)\le 0 \quad \mbox{in}\quad \tilde O.$$
By maximum principle
$$v_i+\phi>(1-3\epsilon_1)|y|^{2-n} \quad \mbox{in}\quad \tilde
O.$$ By (\ref{aug21e2})
$$v_i(y)>(1-4\epsilon_1)|y|^{2-n}\quad \mbox{in}\quad \tilde O.$$
 Since
$h_{\lambda}=\circ(1)|y|^{2-n}$ for $\lambda\in (1/2, 2)$,
(\ref{aug6e2}) is established.

\medskip

Once the moving sphere process can get started at $\lambda_0$, let
$\bar \lambda$ be the critical position where
$w_{\lambda}+h_{\lambda}$ ceases to be positive in
$\Sigma_{\lambda}$. But because of (\ref{aug14e2}),
(\ref{aug14e3}) and (\ref{aug14e4}), the moving sphere process can
reach $\lambda_1$, i.e. $\bar \lambda\ge \lambda_1$. Note that
$T_{\lambda}h_{\lambda}+E_{\lambda}$ only needs to be non positive
in $O_{\lambda}$ because by (\ref{aug14e3})
$w_{\lambda}+h_{\lambda}>0$ in $\Sigma_{\lambda}\setminus
O_{\lambda}$. Since $\lambda_1>1$, by letting $i\to \infty$ we
have
$$U(y)-U^{\lambda_1}(y)\le 0,\quad |y|\ge \lambda_1, $$
which is contradictory to (\ref{aug17e1}).

\begin{rem}
 What is described above is a general procedure of the
application
 of moving sphere method, which will be used a few times in the sequel. Even
 though $\Sigma_{\lambda}$ and $h_{\lambda}$ will be different in different contexts,
 the important thing is to construct
 $h_{\lambda}$ that satisfies (\ref{aug14e2}), (\ref{aug14e3})
 and (\ref{aug14e4}). The way to start the moving sphere process
 and to apply the maximum principle to get a contradiction from
 the standard bubbles are just the same and will not be repeated.
 \end{rem}

To construct $h_{\lambda}$ in this subsection we use the following
crude estimate of $E_{\lambda}$:
\begin{equation}
 \label{aug3e6}
 |E_{\lambda}|\le CM_i^{-\frac 4{n-2}}r^{2-n}+CM_i^{-\frac
 2{n-2}}r^{-n-1}\le CM_i^{-\frac 2{n-2}}r^{2-n},\quad \Sigma_{\lambda}.
\end{equation}
Note that we use $r$ to represent $|y|$. The construction of
$h_{\lambda}$ in this subsection is not subtle with respect to
$\lambda$, we just set $\lambda_0=\frac 12$ and $\lambda_1=2$. We
need the following non positive function: For $2<\alpha<n$, let
$$f_{n,\alpha}(r)=-\frac{1}{(n-\alpha)(2-\alpha)}(r^{2-\alpha}-\lambda^{2-\alpha})
-\frac{\lambda^{n-\alpha}}{(n-\alpha)(n-2)}(r^{2-n}-\lambda^{2-n}).$$
By direct computation one verifies that
\begin{equation}
\label{aug8e4} \left\{\begin{array}{ll} \Delta
f_{n,\alpha}(r)=f_{n,\alpha}''(r)+\frac{n-1}rf_{n,\alpha}'(r)=-r^{-\alpha},
\quad r\ge \lambda,\\
f_{n,\alpha}(\lambda)=f'_{n,\alpha}(\lambda)=0,\\
0\le -f_{n,\alpha}(r)\le C(n,\alpha).
\end{array}
\right.
\end{equation}
This function is mainly used to control minor terms. We define
$h_{\lambda}$ as
$$h_{\lambda}=QM_i^{-\frac 2{n-2}}f_{n,n-2}(\frac{r}{\lambda})$$
where $Q$ is a large number to be determined. Here we see that by
(\ref{aug8e4}), (\ref{aug14e2}) and (\ref{aug14e3}) hold. Also
$h_{\lambda}\le 0$ in $\Sigma_{\lambda}$. Now we verify
(\ref{aug14e4}). First by choosing $Q$ large enough we have
$$\Delta h_{\lambda}+E_{\lambda}\le
-\frac Q2M_i^{-\frac 2{n-2}}r^{2-n},\quad \Sigma_{\lambda}.$$

Since $h_{\lambda}$ is non positive, the term $n^*K_i(M_i^{-\frac
2{n-2}}y+x_i)\xi_{\lambda}^{\frac 4{n-2}}h_{\lambda}$ is also non
positive. The only thing we need to verify is
$$-\mu_i M_i^{-\frac 4{n-2}}h_{\lambda}\le
\frac Q4M_i^{-\frac 2{n-2}}r^{2-n}\quad \mbox{in}\quad
\Sigma_{\lambda}.$$  By direct computation this holds. Proposition
\ref{aug3p1} is established. $\Box$

\bigskip

Next we establish the closeness between $v_i$ and $U_i$, which
satisfies:
$$\Delta U_i+K_i(x_i)U_i^{n^*}=0,\quad \Bbb R^n.\quad U_i(0)=1=\max_{\Bbb R^n} U_i.$$

\begin{prop}
 \label{aug4p2}
 There exist $\delta_1>0$ and $C>0$ such that
 $$|v_i(y)-U_i(y)|\le CM_i^{-\frac{2}{n-2}},\quad
 |y|\le \delta_1 M_i^{\frac{2}{(n-2)^2}}. $$
 \end{prop}

\noindent{\bf Proof of Proposition \ref{aug4p2}}

The proof of Proposition \ref{aug4p2} consists of two steps. First
we show that there exists $\delta_2(n,C_0)>0$ small, so that
\begin{equation}
\label{aug14e6}
v_i(y)\le CU_i(y),\quad |y|\le \delta_2
M_i^{\frac{2}{(n-2)^2}}.
\end{equation}

  The proof of (\ref{aug14e6}) is very similar to that of Lemma 3.2 in
  \cite{CSL2}. The only difference comes from the
  extra term $-\mu_i M_i^{-\frac 4{n-2}}v_i$. For this we let $G_1$
  be the Green's function of the operator
  $-\Delta+\mu_i M_i^{-\frac 4{n-2}}$ with respect to the Dirichlet
  condition on $\Omega_i$ (Recall that $\Omega_i:=B(0,\frac 1{10}M_i^{\frac
  2{n-2}})$). i.e.
  $$\left\{\begin{array}{ll}
  (-\Delta_y +\mu_i M_i^{-\frac 4{n-2}})G_1(x,y)=\delta_x, \quad
  \Omega_i,\\
  G_1(x,y)=0, \quad y\in \partial \Omega_i.
  \end{array}
  \right.
  $$
  By direct computation
  \begin{equation}
  \label{aug14e7}
  G_1(x,y)=\frac 1{\omega_n(n-2)}|x-y|^{2-n}+\phi(x,y)
  \end{equation}
  where $\omega_n$ is the area of $S^{n-1}$, $\phi(x,y)$ satisfies
  \begin{equation}
  \label{aug14e8}
  |\phi(x,y)|\le \mu_i M_i^{-\frac 4{n-2}}C(n)|x-y|^{4-n}.
  \end{equation}
  Once we have this, the rest of the proof is very similar to
  that of lemma 3.2 in \cite{CSL2}.  For $\epsilon>0$ small to be determined, there exists
some constant $\delta_1\in (0,1)$, independent of $i$, such that
for large $i$, let $y_1$ be a minimum of $v_i$ on
$|y|=\delta_1M_i^{\frac 2{(n-2)^2}}$, the following estimates
hold:
\begin{eqnarray*}
v_i(y_1)&\ge &\int_{\Omega_i} G_1(y_1,\eta)K_i(M_i^{-\frac
2{n-2}}y+x_i) v_i(\eta)^{\frac{n+2}{n-2}}d\eta
\\
&\ge&
\int_{B(0,\delta_1M_i^{\frac{2}{(n-2)^2}})}G_1(y_1,\eta)K_i(M_i^{-\frac
2{n-2}}y+x_i) v_i(\eta)^{\frac{n+2}{n-2}}d\eta,
\end{eqnarray*}
 and, using (\ref{aug14e7}) and (\ref{aug14e8}),
$$G_1(y_1,\eta)\ge \frac{(1-\epsilon/2)}{(n-2)\omega_n}|y_1-\eta |^{2-n}
\ge \frac{(1-3\epsilon/4)}{(n-2)\omega_n}|y_1|^{2-n}, \qquad
|\eta|\le \delta_2 |y_1|, $$  if $\delta_2$ is chosen small
enough. Now we use $\lim_{i\to \infty}K_i(x_i)=n(n-2)$ to get
$$v_i(y_1)\ge \frac{(1-\epsilon)n}{\omega_n}|y_1|^{2-n}
\int_{B(0,\delta_2|y_1|)}v_i^{\frac{n+2}{n-2}}(\eta)d\eta $$ On
the other hand, by Proposition \ref{aug3p1}
$$
v_i(y_1)\le  (1+\epsilon) |y_1|^{2-n}.
$$
So
$$\int_{B(0,\delta_2|y_1|)}v_i^{\frac{n+2}{n-2}}(\eta)d\eta\le
(1+4\epsilon)\omega_n/n.$$ A direct computation gives,
$$\int_{\Bbb R^n}U^{\frac{n+2}{n-2}}=\frac{ \omega_n}{ n }.
$$
By the convergence of $v_i$ to $U$, there exists  some $R_1$,
depending only on $n$ and $\epsilon$,
 such that, for large $k$,
$$\int_{R_1\le |\eta|\le \delta_2|y_1|}v_i^{\frac{n+2}{n-2}}d\eta \le \frac{5\epsilon}{n}\omega_n.
$$
Since $v_i\le C_2$ (by(\ref{aug15e1}))
$$\int_{R_1\le |\eta|\le \delta_2|y_1|}v_i^{\frac{2n}{n-2}}d\eta \le
C_2\int_{R_1\le |\eta|\le \delta R_k}v_i^{\frac{n+2}{n-2}}d\eta
\le C\epsilon.$$

For each $2R_1<r<\delta_2|y_1|/2$, we consider $\tilde
v_i(z)=r^{\frac{n-2}2}v_i(rz)$ for $1/2<|z|<2$. Then  $\tilde v_i$
satisfies
$$\Delta \tilde v_i(z)-\mu_i M_i^{-\frac 4{n-2}}r^2\tilde v_i(z)+
K_i(M_i^{-\frac 2{n-2}}rz+x_i)\tilde v_i(z)^{\frac{n+2}{n-2}}
=0,\quad 1/2<|z|<2.$$ We know that $\int_{ \frac 12\le |z|\le 2}
\tilde v_i(z)^{ \frac {2n}{n-2} }\le C\epsilon$. Fix some
universally small $\epsilon>0$, we apply the Moser iteration
technique to obtain $\tilde v_i(z)\le C$ for $\frac 34\le |z|\le
\frac 43$, where $C$ is independent of $k$. With this, we apply
the Harnack inequality to obtain $\max_{ |z|=1} \tilde v_i(z) \le
C\min_{ |z|=1} \tilde v_i(z)$, i.e., $\max_{ |y|=r} v_i(y) \le
C\min_{ |y|=r} v_i(y)$. Then (\ref{aug14e6}) is established.

The second part of the proof is essentially the argument of
 Lemma 3.3 in \cite{CSL2}.
 We state the outline here. Let $w_i=v_i-U_i$. Here we
 recall that $U_i$ satisfies
 $$\left\{\begin{array}{ll}
 \Delta U_i+K_i(x_i)U_i^{n^*}=0 \quad \Bbb R^n,\\
 \\
 U_i(0)=1=\max_{\Bbb R^n}U_i.
 \end{array}
 \right.
 $$
 The equation for $w_i$ is
 \begin{equation}
 \label{aug6e3}
 \left\{\begin{array}{ll}
 (\Delta -\mu_i M_i^{-\frac 4{n-2}}+n^*K_i(M_i^{-\frac
 2{n-2}}y+x_i)\xi_i^{\frac 4{n-2}})w_i=E_i,\quad
 r\le \delta M_i^{\frac{2}{(n-2)^2}}\\
 w_i(0)=0=|\nabla w_i(0)|
 \end{array}
 \right.
 \end{equation}
 where $\xi_i$ is obtained from the mean value theorem:
 $$\xi_i^{\frac 4{n-2}}=\int_0^1(tv_i+(1-t)U_i)^{\frac 4{n-2}}dt$$
 and
 $$E_i=\mu_i M_i^{-\frac 4{n-2}}U_i+
 (K_i(x_i)-K_i(M_i^{-\frac 2{n-2}}y+x_i))U_i^{n^*}.$$
 For $E_i$ we clearly have
 $$|E_i(y)|\le CM_i^{-\frac 4{n-2}}(1+|y|)^{2-n}+CM_i^{-\frac
 2{n-2}}(1+|y|)^{-1-n}.$$
 Let
 $$\Lambda_i=\max_{B(0,\delta M_i^{\frac 2{(n-2)^2}})} \frac{|w_i(y)|}{M_i^{-\frac
 2{n-2}}}.$$
 The goal is to prove $\Lambda_i\le C$. We shall prove by
 contradiction. Suppose $\Lambda_i\to \infty$, let $y_i$ be the
 point that $\lambda_i$ is attained. Let
 $$\bar w_i=\frac{w_i}{\Lambda_iM_i^{-\frac 2{n-2}}}.$$
 Then if $|y_i|$ is bounded, a subsequence of $\bar w_i$ will
 converge to $w$ that satisfies
 \begin{equation}
 \label{aug6e4}
 \left\{\begin{array}{ll}
 \Delta w+n^*U^{\frac 4{n-2}}w=0,\quad \Bbb R^n,\\
 w(0)=0=|\nabla w(0)|, \quad |w(y)|\le 1.
 \end{array}
 \right.
 \end{equation}
 The only function that satisfies (\ref{aug6e4}) is $0$ (Lemma 2.4
 of \cite{CSL2}). This violates $\bar w_i(y_i)=\pm 1$. This
 contradiction forces us to assume $y_i\to \infty$. However by the
 Green's representation theorem, the estimate of $E_i$ makes it
 impossible to have $|\bar w_i(y_i)|=1$. Proposition \ref{aug4p2}
 is established. $\Box$

\subsection{Estimate of $v_i$ over $B(0,\delta M_i^{\frac
4{(n-2)^2}})$}

 Proposition \ref{aug4p2} enables us to improve the estimate of the error term.
  In fact, from the equation for $w_i$ we first have
 \begin{equation}
 \label{aug6e5}
 |w_i(y)|\le CM_i^{-\frac 2{n-2}}|y|^2,\quad |y|\le 10.
 \end{equation}
 Consequently
 $$v_i^{\lambda}(y)=U_i^{\lambda}(y)+O(M_i^{-\frac
 2{n-2}})|y|^{-n}.$$
 So we can write $E_{\lambda}$ as (see (\ref{aug3e4}))
\begin{eqnarray}
E_{\lambda}&=&\mu_i M_i^{-\frac 4{n-2}}U_i^{\lambda}(1-(\frac
{\lambda}r)^4)+(K_i(M_i^{-\frac
2{n-2}}y^{\lambda}+x_i)-K_i(M_i^{-\frac
2{n-2}}y+x_i))(U_i^{\lambda})^{n^*}\nonumber\\
&&+O(M_i^{-\frac 6{n-2}}r^{-n})+O(M_i^{-\frac 4{n-2}}r^{-3-n}).
\label{aug8e1}
\end{eqnarray}
Note that the last two terms come from the difference between
$v_i^{\lambda}$ and $U_i^{\lambda}$. These terms will be estimated
again later. Now by the Taylor expansion of $K_i$ we have
\begin{equation}
\label{aug7e1} E_{\lambda}=E_1+E_2+O(M_i^{-\frac 6{n-2}}r^{1-n})
+O(M_i^{-\frac 4{n-2}}r^{-3-n}).
\end{equation}
where
\begin{equation}
\label{aug9e1} E_1=\mu_i M_i^{-\frac
4{n-2}}U_i^{\lambda}(1-(\frac{\lambda}r)^4) +\frac
1{2n}M_i^{-\frac 4{n-2}}\Delta
K_i(x_i)(\frac{\lambda^4}{r^4}-1)r^2(U_i^{\lambda})^{n^*}
\end{equation}
\begin{eqnarray*}
E_2&=&M_i^{-\frac
2{n-2}}\sum_{j}\partial_jK_i(x_i)(\frac{\lambda^2}{r^2}-1)r\theta_j(U_i^{\lambda})^{n^*}\\
&&+M_i^{-\frac 4{n-2}}\sum_{j\neq
l}\partial_{jl}K_i(x_i)(\frac{\lambda^4}{r^4}-1)r^2\theta_j\theta_l(U_i^{\lambda})^{n^*}\\
&&+\frac 12M_i^{-\frac 4{n-2}}\sum_j\partial_{jj}K_i(x_i)
(\frac{\lambda^4}{r^4}-1)r^2(\theta_j^2-\frac
1n)(U_i^{\lambda})^{n^*}
\end{eqnarray*}
where $\theta_j=y_j/r$. Note that the term of the order
$M_i^{-\frac 6{n-2}}$ has been changed due to the expansion of
$K_i$. Each term in $E_2$ can be considered as a product of a
radial function and an angular function. Each angular function is
an eigenfunction of $-\Delta_{\theta}$ on $S^{n-1}$ ($\theta_j$
corresponds to eigenvalue $n-1$, $\theta_j^2-\frac 1n$ corresponds
to eigenvalue $2n$). By using the ideas in \cite{LiZhang1}
\cite{z1} we construct test functions of the same form. The
current purpose is to prove Proposition \ref{aug7p1} in the
sequel.

Before we state Proposition \ref{aug7p1} we include here a
proposition whose proof can be found in \cite{LiZhang1}:

\begin{prop}
\label{aug8p1} For each $s=1,2$, there exists a unique $C^2$
radial function $g_s$ that satisfies
$$\left\{\begin{array}{ll}
g_s''+\frac{n-1}rg_s'+(n^*K_i(x_i)\tilde\xi_{\lambda}^{\frac
4{n-2}}-\frac{\bar
\lambda_s}{r^2})g_s=r^s((\frac{\lambda}{r})^{2s}-1)(U_i^{\lambda})^{n^*},
\quad \lambda<r<M_i^{\frac 2{n-2}},\\
g_s(\lambda)=0,\quad g_s(M_i^{\frac 2{n-2}})=0.
\end{array}
\right.
$$
where $\lambda\in (1-\epsilon(n),1+\epsilon(n))$, $\epsilon(n)$ is
small, $\bar \lambda_s=s(s+n-2)$, $\tilde \xi_{\lambda}$ is
$$\tilde \xi_{\lambda}^{\frac
4{n-2}}=\int_0^1(tU_i+(1-t)U_i^{\lambda})^{\frac 4{n-2}}dt.$$
Moreover, there exists a dimensional constant $C_0(n)$ so that
\begin{equation}
\label{aug8e2} 0\le g_s(r)\le C_0(1-\frac{\lambda}r)r^{2-n}, \quad
\lambda<r<M_i^{\frac 2{n-2}}.
\end{equation}
\end{prop}

By comparing $\tilde \xi_{\lambda}$ and $\xi_{\lambda}$ we see
that $\tilde \xi_{\lambda}$ is radial and is very close to
$\xi_{\lambda}$ for $r\le \delta_2 M_i^{\frac 2{(n-2)^2}}$. For
$r\ge \delta_2 M_i^{\frac 2{(n-2)^2}}$ both terms are comparable
to $r^{2-n}$.

 Next we show
\begin{prop}
\label{aug7p1} Given $\epsilon>0$, there exists $\delta_3>0$ such
that for all large $i$,
$$\min_{|y|=r}v_i(y)\le (1+\epsilon)r^{2-n},\quad
10\le r\le \delta_3 M_i^{\frac{4}{(n-2)^2}}.$$
\end{prop}

\noindent{\bf Proof of Proposition \ref{aug7p1}:} We prove this by
a contradiction. Suppose there exist $\epsilon_0>0$ and
$\epsilon_i\to 0$ such that
\begin{equation}
\label{aug7e2} \min_{|y|=r_i}v_i(y)\ge (1+\epsilon_0)r_i^{2-n},
\quad \mbox{for some } r_i\in (\delta_2M_i^{\frac{2}{(n-2)^2}},
\epsilon_i M_i^{\frac{4}{(n-2)^2}}).
\end{equation}
By the convergence of $v_i$ to $U$, $r_i\to \infty$. Let
$\Sigma_{\lambda}=B(0,\epsilon_iM_i^{\frac {4}{(n-2)^2}})\setminus
\bar B_{\lambda}$
 and
 \begin{equation}
 \label{aug9e7}
h_1=-M_i^{-\frac 2{n-2}}\sum_j\partial_jK_i(x_i)\theta_jg_1(r).
\end{equation}
\begin{equation}
 \label{aug9e8}
h_2=-M_i^{-\frac 4{n-2}}g_2(r)(\sum_{j\neq
l}\partial_{jl}K_i(x_i)\theta_j\theta_l+\frac
12\sum_j\partial_{jj}K_i(x_i)(\theta_j^2-\frac 1n)).
\end{equation}
Then we have
\begin{equation}
\label{aug8e3} (\Delta +n^*K_i(x_i)\tilde \xi_{\lambda}^{\frac
4{n-2}}) (h_1+h_2)+E_2=0
\end{equation}
and, by Proposition \ref{aug8p1}
\begin{eqnarray}
|h_1(y)|&\le & C|\nabla K_i(x_i)|M_i^{-\frac
2{n-2}}r^{2-n},\nonumber\\
|h_2(y)|&\le &CM_i^{-\frac 4{n-2}}r^{2-n}. \label{aug17e2}
\end{eqnarray}

Let
$$h_3(r)=QM_i^{-\frac 4{n-2}}f_{n,3}$$
where $Q>>1$ is to be determined. Then first we notice that
$h_3<0$ in $\Sigma_{\lambda}$. Also by the definition of
$f_{n,\alpha}$ we have
\begin{equation}
\label{aug17e3}
 \Delta h_3=-QM_i^{-\frac 4{n-2}}r^{-3}, \quad
\Sigma_{\lambda}.
\end{equation}

By (\ref{aug17e2}) and (\ref{aug8e4}), each of $h_j, j=1,2,3$
satisfies (\ref{aug14e2}) and (\ref{aug14e3}). Now by
(\ref{aug8e3}) and (\ref{aug17e3})
\begin{eqnarray}
&&T_{\lambda}(h_1+h_2+h_3)\nonumber\\
&=&-E_2-\mu_i M_i^{-\frac 4{n-2}}(h_1+h_2+h_3)\nonumber\\
&&+n^*(K_i(M_i^{-\frac 2{n-2}}y+x_i)\xi_{\lambda}^{\frac
4{n-2}}-K_i(x_i)\tilde \xi_{\lambda}^{\frac 4{n-2}})(h_1+h_2)\nonumber\\
&&-QM_i^{-\frac 4{n-2}}r^{-3}+n^*K_i(M_i^{-\frac
2{n-2}}y+x_i)\xi_{\lambda}^{\frac 4{n-2}}h_3. \label{aug8e5}
\end{eqnarray}
Note that since $h_3\le 0$ the last two terms have a good sign. To
estimate other terms in (\ref{aug8e5}) we first use (\ref{aug8e2})
and (\ref{aug8e4}) to get
\begin{equation}
\label{aug8e6} -\mu_i M_i^{-\frac
4{n-2}}(h_1+h_2+h_3)=O(M_i^{-\frac 6{n-2}}r^{2-n})+O(M_i^{-\frac
8{n-2}}).
\end{equation}

 Next by Proposition \ref{aug4p2} we estimate the following:
 \begin{eqnarray}
 &&K_i(M_i^{-\frac 2{n-2}}y+x_i)\xi_{\lambda}^{\frac
 4{n-2}}-K_i(x_i)\tilde \xi_{\lambda}^{\frac 4{n-2}}\nonumber\\
 &=&(K_i(M_i^{-\frac 2{n-2}}y+x_i)-K_i(x_i))\xi_{\lambda}^{\frac
 4{n-2}}+K_i(x_i)(\xi_{\lambda}^{\frac 4{n-2}}-\tilde
 \xi_{\lambda}^{\frac{4}{n-2}})\nonumber\\
 &=&O(M_i^{-\frac{2}{n-2}}r^{-3})+K_i(x_i)(\xi_{\lambda}^{\frac 4{n-2}}-\tilde
 \xi_{\lambda}^{\frac{4}{n-2}}).
 \label{aug9e9}
 \end{eqnarray}
 To estimate the last term of the above, we use
\begin{eqnarray}
\xi_{\lambda}^{\frac{4}{n-2}}&=&\int_0^1(tv_i+(1-t)v_i^{\lambda})^{\frac
4{n-2}}dt \nonumber\\
&=&\int_0^1(tU_i+(1-t)U_i^{\lambda}+a)^{\frac 4{n-2}}dt \nonumber\\
&=&\tilde \xi_{\lambda}^{\frac 4{n-2}}+ar^{n-6}
 \label{aug9e10}
\end{eqnarray}
where
\begin{equation}
\label{aug9e11}
|a|\le \left\{\begin{array}{ll}
CM_i^{-\frac{2}{n-2}},
\quad \lambda<r<\delta M_i^{\frac{2}{(n-2)^2}},\\
Cr^{2-n},\quad r\ge \delta M_i^{\frac{2}{(n-2)^2}}, \quad y\in
O_{\lambda}.
\end{array}
\right.
\end{equation}
Putting (\ref{aug17e2}) (\ref{aug9e9}) (\ref{aug9e10}) and
(\ref{aug9e11}) together we have
\begin{equation}
\label{aug8e7}  n^*(K_i(M_i^{-\frac
2{n-2}}y+x_i)\xi_{\lambda}^{\frac 4{n-2}}-K_i(x_i)\tilde
\xi_{\lambda}^{\frac 4{n-2}})(h_1+h_2) =O(M_i^{-\frac
4{n-2}}r^{-3}), \quad O_{\lambda}.
\end{equation}
Thus, by (\ref{aug3e5}), (\ref{aug7e1}), (\ref{aug8e5}),
(\ref{aug8e6}) and (\ref{aug8e7}) we have
\begin{eqnarray*}
&&T_{\lambda}(w_{\lambda}+h_1+h_2+h_3) \\
&\le &E_1+O(M_i^{-\frac 6{n-2}}r^{2-n})+O(M_i^{-\frac
4{n-2}}r^{-3})\\
&&+O(M_i^{-\frac 8{n-2}})-QM_i^{-\frac 4{n-2}}r^{-3},\quad
\mbox{in}\quad O_{\lambda}.
\end{eqnarray*}
It is easy to verify that
$$|E_1(y)|\le CM_i^{-\frac 4{n-2}}r^{2-n}, \quad M_i^{-\frac
8{n-2}}=\circ(1)M_i^{-\frac 4{n-2}}r^{-3}.$$
So by choosing $Q$
large enough we have
$$T_{\lambda}(w_{\lambda}+h_1+h_2+h_3)\le 0 \quad O_{\lambda}.$$
Proposition \ref{aug7p1} is established.

\begin{rem} In the proof of Proposition \ref{aug7p1} we don't need
the sign of $\Delta K_i(x_i)$.
\end{rem}

\subsection{The vanishing rate of $|\nabla K_i(x_i)|$} Next we
improve the estimate of $v_i-U_i$: First we have
\begin{prop}
\label{aug9p1} There exist $\delta_4>0$ and $C>0$ such that
$$v_i(y)\le CU_i(y)\quad |y|\le \delta_4 M_i^{\frac{4}{(n-2)^2}}.$$
\end{prop}
The proof is similar to the step one of Proposition \ref{aug4p2}.
$\Box$
\bigskip

 To
further estimate $v_i-U_i$ more precisely, we need the following
Pohozaev Identity for
$$\Delta f(x)-t f(x)+H(x)f(x)^{n^*}=0\quad B_{\sigma}.$$
\begin{eqnarray*}
&&\int_{B_{\sigma}}(\frac{n-2}{2n}(\nabla H\cdot x)f^{\frac{2n}{n-2}}-t f^2)\\
&=&\int_{\partial B_{\sigma}}(\frac{n-2}{2n}\sigma
Hf^{\frac{2n}{n-2}}+\sigma |\frac{\partial f}{\partial
\nu}|^2-\frac{\sigma}2|\nabla f|^2+\frac{n-2}2\frac{\partial
f}{\partial \nu}f-\frac{t}2\sigma f^2).
\end{eqnarray*}

Let $\tilde v_i(y)=v_i(y+e)$, where $e=\frac{\nabla
K_i(x_i)}{|\nabla K_i(x_i)|}$ is a unit vector. Let
$$\tilde K_i(y)=K_i(M_i^{-\frac{2}{n-2}}(y+e)+x_i),$$ then $\tilde v_i$
satisfies
$$\Delta \tilde v_i(y)-\mu_i M_i^{-\frac 4{n-2}}\tilde v_i(y)+\tilde K_i(y)\tilde v_i(y)^{n^*}=0,\quad
|y|\le L_i:=\frac{\delta_4}{2}M_i^{\frac{4}{(n-2)^2}}.$$
So the
Pohozaev Identity applied to $\tilde v_i$ over $B_{L_i}$ gives
\begin{eqnarray}
&&\int_{B_{L_i}}\bigg (\frac{n-2}{2n}(\nabla \tilde K_i(y)\cdot y)
\tilde v_i^{\frac{2n}{n-2}}-\mu_i M_i^{-\frac 4{n-2}}\tilde v_i^2(y)\bigg )dy \nonumber\\
&=&\int_{\partial B_{L_i}}(\frac{n-2}{2n}L_i \tilde K_i(y)\tilde
v_i^{\frac{2n}{n-2}} +L_i |\frac{\partial \tilde v_i}{\partial
\nu}|^2-\frac{L_i}2|\nabla \tilde v_i|^2\nonumber \\
&&+\frac{n-2}2\frac{\partial \tilde v_i}{\partial \nu}\tilde
v_i-\frac{\mu_i}2M_i^{-\frac 4{n-2}}L_i\tilde v_i^2). \label{poho}
\end{eqnarray}
By Proposition \ref{aug9p1} and standard elliptic estimates
 the right hand side of (\ref{poho}) is
$O(M_i^{-\frac{4}{n-2}})$. Then by using $e=\nabla
K_i(x_i)/|\nabla K_i(x_i)|$ we see that the left hand side of the
Pohozaev identity is greater than
$$C|\nabla K_i(x_i)|M_i^{-\frac{2}{n-2}}+O(M_i^{-\frac 4{n-2}})$$
for some $C>0$. Consequently
\begin{equation}
|\nabla K_i(x_i)|\le CM_i^{-\frac{2}{n-2}}. \label{mar29e1}
\end{equation}
Base on (\ref{mar29e1}) we can write the equation for
$w_i:=v_i-U_i$ as
$$(\Delta +n^*K_i(M_i^{-\frac 2{n-2}}y+x_i)\xi_i^{\frac
4{n-2}})w_i=O(M_i^{-\frac {4}{n-2}})(1+r)^{2-n},\quad r\le
\delta_4 M_i^{\frac 4{(n-2)^2}}.$$ Then the same estimate in
Proposition \ref{aug4p2} gives
\begin{equation}
\label{aug9e5} |v_i(y)-U_i(y)|\le CM_i^{-\frac 4{n-2}},\quad
|y|\le \delta_4 M_i^{\frac 4{(n-2)^2}}.
\end{equation}
Using the fact that $v_i(0)=U_i(0)$ and $\nabla v_i(0)=\nabla
U_i(0)=0$, we have
\begin{equation}
\label{aug14e18} |v_i(y)-U_i(y)|\le CM_i^{-\frac 4{n-2}}|y|^2,
\quad |y|\le 10.
\end{equation}

\subsection{Harnack inequality on $B(0,\delta M_i^{\frac 2{n-2}})$}

 Now we establish
\begin{prop}
\label{aug9p2}  For any $\epsilon>0$, there exists $\delta_5>0$
depending on $\epsilon$ and $n$ such that
$$\min_{|y|=r}v_i(y)\le (1+\epsilon)r^{2-n},\quad
r\le \delta_5 M_i^{\frac 2{n-2}}. $$
\end{prop}

\noindent{\bf Proof of Proposition \ref{aug9p2}:} We still prove
it by a contradiction by assuming that there exist $\epsilon_0$
and $\epsilon_i\to 0$ such that
\begin{equation}
\label{aug9e12}
\min_{|y|=\epsilon_iM_i^{\frac 2{n-2}}} v_i(y) \ge
(1+\epsilon_0)|y|^{2-n}.
\end{equation}

 As a consequence of (\ref{aug14e18}) (see also (\ref{aug6e5}))
$$v_i^{\lambda}(y)=U_i^{\lambda}(y)+O(M_i^{-\frac
4{n-2}}|y|^{-n}).$$

Therefore, in stead of (\ref{aug7e1}) we now have
\begin{equation}
\label{aug9e6} E_{\lambda}=E_1+E_2+O(M_i^{-\frac
{8}{n-2}}r^{1-n})(1-\frac{\lambda}r) +O(M_i^{-\frac
6{n-2}}r^{-3-n})(1-\frac{\lambda}r).
\end{equation}

Note that we include $1-\frac{\lambda}r$ deliberately because the
dominant term now vanishes on $\partial B_{\lambda}$. We still
construct $h_1$ and $h_2$ as in (\ref{aug9e7}) and (\ref{aug9e8}).
Because of the new rate of $|\nabla K_i(x_i)|$ we now have
$$|h_1|+|h_2|=O(M_i^{-\frac{4}{n-2}}r^{2-n}).$$
We note that (\ref{aug8e3}) also holds. Now we have
\begin{eqnarray*}
T_{\lambda}(h_1+h_2)&=&-E_2-\mu_i M_i^{-\frac 4{n-2}}(h_1+h_2)\\
&&+n^*(K_i(M_i^{-\frac 2{n-2}}y+x_i)\xi_{\lambda}^{\frac
4{n-2}}-K_i(x_i)\tilde \xi_{\lambda}^{\frac 4{n-2}})(h_1+h_2).
\end{eqnarray*}
By the new rate of $|\nabla K_i(x_i)|$ we have
$$|\mu_i M_i^{-\frac{4}{n-2}}(h_1+h_2)|\le CM_i^{-\frac
8{n-2}}r^{2-n}(1-\frac{\lambda}r).$$

Also by using this new rate of $|\nabla K_i(x_i)|$ in
(\ref{aug9e9}) we have
\begin{eqnarray*}
&&K_i(M_i^{-\frac 2{n-2}}y+x_i)\xi_{\lambda}^{\frac 4{n-2}}
-K_i(x_i)\tilde \xi_{\lambda}^{\frac 4{n-2}}\\
&=&O(M_i^{-\frac 4{n-2}}r^{-2})+K_i(x_i)(\xi_{\lambda}^{\frac
4{n-2}}-\tilde \xi_{\lambda}^{\frac 4{n-2}}).
\end{eqnarray*}
Corresponding to (\ref{aug9e10}) and (\ref{aug9e11}) we now have
$$\xi_{\lambda}^{\frac{4}{n-2}}=\tilde
\xi_{\lambda}^{\frac{4}{n-2}}+ar^{n-6} \quad O_{\lambda}$$ where
$$|a|\le \left\{\begin{array}{ll}
O(M_i^{-\frac 4{n-2}}),\quad \lambda<r<\delta
M_i^{\frac{4}{(n-2)^2}},\\
Cr^{2-n},\quad \delta M_i^{\frac{4}{(n-2)^2}}\le r\le
\epsilon_iM_i^{\frac 2{n-2}}.
\end{array}
\right.
$$
Consequently
\begin{eqnarray*}
&&n^*(K_i(M_i^{-\frac 2{n-2}}y+x_i)\xi_{\lambda}^{\frac
4{n-2}}-K_i(x_i)\tilde \xi_{\lambda}^{\frac 4{n-2}})(h_1+h_2)\\
&=&O(M_i^{-\frac 8{n-2}}r^{2-n})(1-\frac{\lambda}r)+\tilde E_3.
\end{eqnarray*}
where $\tilde E_3$ is $0$ when $r\le \delta_4 M_i^{\frac
4{(n-2)^2}}$ and is $O(M_i^{-\frac 4{n-2}}r^{-2-n})$ when $r\ge
\delta_4 M_i^{\frac 4{(n-2)^2}}$.

Now we can combine $w_{\lambda}$ with $h_1, h_2$:
\begin{eqnarray*}
&&T_{\lambda}(w_{\lambda}+h_1+h_2)\\
&\le &E_1+O(M_i^{-\frac 8{n-2}}r^{2-n})(1-\frac{\lambda}r)
+O(M_i^{-\frac 6{n-2}}r^{-3-n})(1-\frac{\lambda}r)+\tilde E_3.
\end{eqnarray*}
Now the second term in $E_1$ becomes the dominant term. In fact,
since $\Delta K_i(x_i)$ is large we have
\begin{eqnarray}
&&T_{\lambda}(w_{\lambda}+h_1+h_2)\nonumber \\
&\le &\mu_i M_i^{-\frac
4{n-2}}U_i^{\lambda}(1-(\frac{\lambda}r)^4) -\frac
1{3n}M_i^{-\frac 4{n-2}}\Delta
K_i(x_i)(1-(\frac{\lambda}r)^4)r^2(U_i^{\lambda})^{n^*}.
\label{aug9e13}
\end{eqnarray}

Note that $\Delta K_i(x_i)$ makes the right hand side of
(\ref{aug9e13}) negative when $r$ is close to, or comparable to
$\lambda$. But the right hand side becomes positive when $r$ is
large. So we need to construct the following function to deal with
this. Let
$$h_3(y)=\int_{\Sigma_{\lambda}}\mu_i M_i^{-\frac
4{n-2}}G(y,\eta)U_i^{\lambda}(\eta)(1-(\frac {\lambda}{|\eta
|})^4)d\eta .
$$

For $h_3$ we have
$$0\le h_3(y)\le C(C_0, n)M_i^{-\frac 4{n-2}}|y|^{4-n}, $$ (so $h_3=\circ(1)|y|^{2-n}$)
and
$$-\Delta h_3=\mu_i M_i^{-\frac 4{n-2}}U_i^{\lambda}(1-(\frac {\lambda}{|\eta
|})^4),\quad \Sigma_{\lambda}.
$$
Because of this equation, the bad term now becomes
$n^*K_i(M_i^{-\frac 2{n-2}}y+x_i)\xi_{\lambda}^{\frac 4{n-2}}h_3$.
To be more precise, we have
\begin{eqnarray*}
&&T_{\lambda}(w_{\lambda}+\sum_{s=1}^3h_s)\\
&\le &-\frac 1{3n}M_i^{-\frac 4{n-2}}\Delta K_i(x_i)
(1-(\frac{\lambda}{r})^4)r^2(U_i^{\lambda})^{n^*}+n^*K_i(M_i^{-\frac
2{n-2}}y+x_i)\xi_{\lambda}^{\frac 4{n-2}}h_3\\
&\le &-\frac 1{3n}M_i^{-\frac 4{n-2}}\Delta K_i(x_i)
(1-(\frac{\lambda}{r})^4)r^2(U_i^{\lambda})^{n^*}\\
&&+C(n,C_0)M_i^{-\frac 4{n-2}}|y|^{4-n}(1+|y|)^{-4}(1-\lambda/r)\\
&\le &0
\end{eqnarray*}
where in the last step we used the largeness of $\Delta K_i(x_i)$.
With this inequality the moving sphere argument applies as before
to get a contradiction. Proposition \ref{aug9p2} is established.
$\Box$

\medskip

\subsection{The completion of the proof of Theorem \ref{thm1}}
 By Proposition \ref{aug9p2} and the step one of the proof of
 Proposition \ref{aug4p2} we have, for some $\delta_5>0$ small and
 $C$ large
 \begin{equation}
 \label{aug15e2}
 v_i(y)\le CU_i(y),\quad |y|\le \delta_5 M_i^{\frac 2{n-2}}.
 \end{equation}

 Note that we can not get better estimate on $v_i-U_i$ as before
 because terms of order $O(M_i^{-\frac 4{n-2}})$ prevent us from
 getting estimates better than (\ref{aug9e5}) and (\ref{aug14e18}).
 By using the vanishing rate of $|\nabla K_i(x_i)|$ (\ref{mar29e1})
 and (\ref{aug15e2}), (\ref{aug9e5}),(\ref{aug14e18}) we shall get
 a contradiction to (\ref{aug14e9}) from the Pohozaev identity as follows.

Let $L_i=\delta_5M_i^{\frac 2{n-2}}$, we apply the Pohozaev
identity to $v_i$ on $B_{L_i}$. Then the right hand side is
$$\int_{\partial B_{L_i}}\{\frac{n-2}{2n}L_iK_i
v_i^{\frac {2n}{n-2}}+L_i|\frac{\partial v_i}{\partial
\nu}|^2-\frac{L_i}2 |\nabla v_i|^2+\frac{n-2}2\frac{\partial
v_i}{\partial \nu}v_i-\frac{\mu_i M_i^{-\frac
4{n-2}}}2L_iv_i^2\}.$$

By (\ref{aug15e2}) and the corresponding gradient estimate the
right hand side of the Pohozaev identity is $O(M_i^{-2})$.

The left hand side of the Pohozaev identity is
$$
\int_{B_{L_i}}\{\frac{n-2}{2n}M_i^{-\frac 2{n-2}}(\nabla
K_i(M_i^{-\frac 2{n-2}}y+x_i)\cdot y)v_i^{\frac {2n}{n-2}}-\mu_i
M_i^{-\frac 4{n-2}}v_i^2\}
$$

We call the first term of the above $L_1$, the second term $L_2$.
Clearly by the convergence of $v_i$ to $U$, we have
$$L_2\ge -\mu_i M_i^{-\frac 4{n-2}}C(n).$$
Now we estimate $L_1$, for which we first have
$$\nabla K_i(M_i^{-\frac 2{n-2}}y+x_i)\cdot y
=\sum_j\partial_jK_i(x_i)y_j+M_i^{-\frac
2{n-2}}\sum_{l,j}\partial_{jl}K_i(x_i)y_jy_l+O(M_i^{-\frac
4{n-2}}r^3).$$

Then we can write $L_1$ as
\begin{eqnarray*}
L_1&=&\frac{n-2}{2n}M_i^{-\frac
2{n-2}}\int_{B_{L_i}}\sum_j\partial_jK_i(x_i)y_jv_i^{\frac{2n}{n-2}}\\
 &&+\frac{n-2}{2n}M_i^{-\frac
4{n-2}}\int_{B_{L_i}}\sum_{jl}\partial_{jl}K_i(x_i)y_jy_lv_i^{\frac{2n}{n-2}}dy
+O(M_i^{-\frac 6{n-2}})\\
&=&L_{11}+L_{12}+O(M_i^{-\frac 6{n-2}}).
\end{eqnarray*}

To estimate $L_{11}$ let $\bar L_i=M_i^{\frac 1{n-2}}$, then
\begin{eqnarray*}
&&\int_{B_{L_i}}y_jv_i^{\frac {2n}{n-2}}=\int_{B_{\bar
L_i}}+\int_{B_{L_i}\setminus B_{\bar L_i}}\\
&=&\int_{B_{\bar L_i}}y_j(U_i+O(M_i^{-\frac 4{n-2}}))^{\frac
{2n}{n-2}}dy+O(M_i^{-1})\\
&=&O(M_i^{-\frac 4{n-2}})+O(M_i^{-1}).
\end{eqnarray*}

So by (\ref{mar29e1}) $L_{11}=\circ (M_i^{-\frac 4{n-2}}).$ By
similar estimate we see that the leading term in $L_{12}$ is
$\frac{n-2}{2n^2}M_i^{-\frac 4{n-2}}\Delta K_i(x_i)\int_{\Bbb
R^n}|y|^2U_i^{\frac{2n}{n-2}}dy$. Then the largeness of $\Delta
K_i(x_i)$ clearly leads to a contradiction in the Pohozaev
identity. Theorem \ref{thm1} is established. $\Box$

\subsection{The Proof of Corollary \ref{cor1}}
Let
$$u(y)=(\frac{2}{1+|y|^2})^{\frac{n-2}2}w(\pi^{-1}(y)),\quad y\in
\Bbb R^n$$ and $g_0$ denote the standard metric on $\Bbb S^n$. In
stereographic projection
$$g_0=\sum_{i=1}^{n+1}dx_i^2=\bigg \{(\frac
2{1+|y|^2})^{\frac{n-2}2}\bigg \}^{\frac 4{n-2}}dy^2. $$ Since the
north pole is not a critical point of $R$, we know $u\sim
O(|y|^{2-n})$ at infinity.
 Then the equation for $u$ becomes
 \begin{equation}
 \label{8feb26e1}
 \left\{\begin{array}{ll}
 \Delta u(y)-\mu(y)(\frac
 2{1+|y|^2})^2u(y)+K(y)u(y)^{\frac{n+2}{n-2}}=0,\quad \Bbb R^n,\\
 u(y)\sim O(|y|^{2-n})\quad \mbox{at }\infty.
 \end{array}
 \right.
 \end{equation}
By applying Theorem \ref{thm1} we obtain
\begin{equation}
\label{8mar4e1} u(y)\le C, \quad y\in \Bbb R^n.
\end{equation}
If $\mu\equiv 0$, the result in \cite{z1} yields (\ref{8mar4e1})
only under the assumption $\Delta K(y)>0$ for each critical point
$y$. The upper bound on $u$ gives the upper bound on $w$, then by
Harnack inequality
$$\frac 1C\le w(x)\le C,\quad \mbox{on }\quad \Bbb S^n. $$
Corollary \ref{cor1} is established. $\Box$

\end{document}